\documentclass[12pt,reqno]{amsart}

\usepackage{amssymb}

\usepackage{mathrsfs}

\DeclareMathAlphabet{\mathpzc}{OT1}{pzc}{m}{it}

\usepackage[all]{xy}

\entrymodifiers={+!!<0pt,\fontdimen22\textfont2>}

\def\sS{\mbox{\sf S}}
\def\sT{\mbox{\sf T}}
\def\sU{\mbox{\sf U}}

\def\adots{\mathinner{\mkern1mu\raise1.0pt\vbox{\kern7.0pt\hbox{.}}\mkern2mu\raise4.0pt\hbox{.}\mkern2mu\raise7.0pt\hbox{.}\mkern1mu}}

\def\EssIm{\operatorname{Ess.\!Im}}

\def\Hom{\operatorname{Hom}}

\def\Ker{\operatorname{Ker}}

\def\prod{\operatorname{prod}}

\numberwithin{equation}{part}

\newtheorem{Lemma}{Lemma}[section]

\newtheorem{Theorem}[Lemma]{Theorem}
\newtheorem{Proposition}[Lemma]{Proposition}

\theoremstyle{definition}

\newtheorem{Remark}[Lemma]{Remark}

\begin{document}

\setlength{\parindent}{0pt}
\setlength{\parskip}{7pt}

\title[Reflecting recollements]
{Reflecting recollements}

\author{Peter J\o rgensen}
\address{School of Mathematics and Statistics,
Newcastle University, Newcastle upon Tyne NE1 7RU,
United Kingdom}
\email{peter.jorgensen@ncl.ac.uk}
\urladdr{http://www.staff.ncl.ac.uk/peter.jorgensen}





\begin{abstract} 
  
  A recollement describes one triangulated category $\sT$ as ``glued
  together'' from two others, $\sS$ and $\sU$.  The definition is not
  symmetrical in $\sS$ and $\sU$, but this note shows how $\sS$ and
  $\sU$ can be interchanged when $\sT$ has a Serre functor.

\end{abstract}

\maketitle


A recollement of triangulated categories $\sS$, $\sT$, $\sU$ is a
diagram of triangulated functors
\begin{equation}
\label{equ:recollement}
\vcenter{
\xymatrix
{
&&\sS
    \ar[rr]^{i_*} & & 
  \sT 
    \ar[rr]^{j^*}
    \ar@/_1.5pc/[ll]_{i^*} \ar@/^1.5pc/[ll]^{i^!} & &
  \sU
    \ar@/_1.5pc/[ll]_{j_!} \ar@/^1.5pc/[ll]^{j_*} & &
}
        }
\end{equation}
satisfying a number of conditions given in Remark
\ref{rmk:recollements} below.

Recollements are important in algebraic geometry and representation
theory, see for instance \cite{BBD}, \cite{Konig}, \cite{Miyachi}.
They were introduced and de\-ve\-lo\-ped in \cite{BBD}, and as
indicated by the terminology, one thinks of $\sT$ as being ``glued
together'' from $\sS$ and $\sU$.  Indeed, in the canonical example of
a recollement, $\sT$ is a derived category of sheaves on a space, and
$\sS$ and $\sU$ are derived categories of sheaves on a closed subset
and its open complement, respectively.  Other examples of a more
algebraic nature can be found in \cite{Konig}.

The recollement \eqref{equ:recollement} is not symmetrical in $\sS$
and $\sU$: There are only two arrows pointing to the right, but four
pointing to the left.  So there is no particular reason to think that
it should be possible to interchange $\sS$ and $\sU$, that is, use
\eqref{equ:recollement} to construct another recollement of the form
\begin{equation}
\label{equ:reflected_recollement}
\vcenter{
\xymatrix
{
&&\sU
    \ar[rr] & & 
  \sT 
    \ar[rr]
    \ar@/_1.5pc/[ll]_{} \ar@/^1.5pc/[ll]^{} & &
  \sS.
    \ar@/_1.5pc/[ll] \ar@/^1.5pc/[ll] & &
}
        }
\end{equation}
Nevertheless, that is precisely what this note does in Theorem
\ref{thm:main} below, under the assumption that $\sT$ has a Serre
functor; see Remark \ref{rmk:Serre} for the definition.

In fact, it will be showed that there are two different ways to get
re\-col\-le\-ments of the form \eqref{equ:reflected_recollement}, one
involving the four upper functors from \eqref{equ:recollement} and
another involving the four lower functors.

For the rest of the note, $k$ will denote a field, and the category
$\sT$ of the recollement \eqref{equ:recollement} will be assumed to be
a skeletally small $k$-linear triangulated category with finite
dimensional $\Hom$-sets and split idempotents.  It will also be
assumed that $\sT$ has a Serre functor $T$.  By $\widetilde{T}$ is
denoted a quasi-inverse to $T$.

Let me start with two remarks explaining the formalism of recollements
and Serre functors.

\begin{Remark}
[Recollements, cf.\ {\cite[sec.\ 1.4]{BBD}}]
\label{rmk:recollements}
The recollement \eqref{equ:recollement} is defined by the following
properties.
\begin{enumerate}

  \item  $(i^*,i_*)$, $(i_*,i^!)$, $(j_!,j^*)$, and $(j^*,j_*)$ are
         pairs of adjoint functors.

\smallskip

  \item  $j^* i_* = 0$.

\smallskip

  \item  $i_*$, $j_!$, and $j_*$ are fully faithful.

\smallskip

  \item  Each object $X$ in $\sT$ determines distinguished triangles
\smallskip
\begin{enumerate}

  \item  $i_*i^!X \longrightarrow X \longrightarrow j_*j^*X
         \longrightarrow \;\;$ and

  \item  $j_!j^*X \longrightarrow X \longrightarrow i_*i^*X
         \longrightarrow$ 

\end{enumerate}
\smallskip
where the arrows into and out of $X$ are counit and unit morphisms of
the relevant adjunctions.
\end{enumerate}
\end{Remark}

\begin{Remark}
[Serre functors, cf.\ {\cite[sec.\ I.1]{RVdB}}]
\label{rmk:Serre}
Let $(-)^{\vee}$ denote the functor $\Hom_k(-,k)$.  A right Serre
functor $F$ for $\sT$ is an endofunctor for which there are natural
isomorphisms 
\[
  \sT(X,Y) \cong \sT(Y,FX)^{\vee},
\]
and a left Serre functor $\widetilde{F}$ is an endofunctor for which
there are natural isomorphisms 
\[
  \sT(X,Y) \cong \sT(\widetilde{F}Y,X)^{\vee}.
\]
A Serre functor is an essentially surjective right Serre functor.

A right Serre functor is fully faithful, and hence a Serre functor is
an autoequivalence.

If there is a right Serre functor $F$ and a left Serre functor
$\widetilde{F}$, then $F$ is in fact a Serre functor and
$\widetilde{F}$ is a quasi-inverse of $F$.
\end{Remark}

On this basis, it is possible to prove that the categories $\sS$ and
$\sU$ in \eqref{equ:recollement} can be interchanged.  First, however,
two propositions which may be of independent interest.

\begin{Proposition}
\label{pro:Serre}
The category $\sS$ has a Serre functor $S$ with quasi-inverse
$\widetilde{S}$:
\[
  S = i^!Ti_* \;\;\;\; \mbox{and} \;\;\;\; \widetilde{S} = i^*\widetilde{T}i_*.
\]
The category $\sU$ has a Serre functor $U$ with quasi-inverse
$\widetilde{U}$:
\[
  U = j^*Tj_! \;\;\;\; \mbox{and} \;\;\;\; \widetilde{U} = j^*\widetilde{T}j_*.
\]
\end{Proposition}

\begin{proof}
By Remark \ref{rmk:Serre}, it is enough to show that $S$ and
$\widetilde{S}$ are, respectively, a right and a left Serre functor
for $\sS$, and similarly for $U$ and $\widetilde{U}$.  This can be
done directly, 
\renewcommand{\arraystretch}{1.2}
\[
  \begin{array}{rcll}
    \sS(Y,SX)^{\vee} 
    & =     & \sS(Y,i^!Ti_*X)^{\vee} 
            & \mbox{by definition} \\
    & \cong & \sT(i_*Y,Ti_*X)^{\vee}
            & \mbox{$i_*$ left-adjoint of $i^!$} \\
    & \cong & \sT(i_*X,i_*Y)
            & \mbox{$T$ right Serre functor} \\
    & \cong & \sS(X,Y)
            & \mbox{$i_*$ fully faithful}
  \end{array}
\]
\renewcommand{\arraystretch}{1}
$\!\!$and
\renewcommand{\arraystretch}{1.2}
\[
  \begin{array}{rcll}
    \sS(\widetilde{S}Y,X)^{\vee}
    & =     & \sS(i^*\widetilde{T}i_*Y,X)^{\vee}
            & \mbox{by definition} \\
    & \cong & \sT(\widetilde{T}i_*Y,i_*X)^{\vee}
            & \mbox{$i_*$ right-adjoint of $i^*$} \\
    & \cong & \sT(i_*X,i_*Y)
            & \mbox{$\widetilde{T}$ left Serre functor} \\
    & \cong & \sS(X,Y)
            & \mbox{$i_*$ fully faithful.}
  \end{array}
\]
\renewcommand{\arraystretch}{1}
$\!\!$Similar computations work for $U$ and $\widetilde{U}$.
\end{proof}

\begin{Proposition}
\label{pro:adjoints}
The functors $i^*$ and $j_!$ have left-adjoint functors given by
\[
  i_! = \widetilde{T}i_*S = \widetilde{T}i_*i^!Ti_* 
  \;\;\;\; \mbox{and} \;\;\;\;
  j^? = \widetilde{U}j^*T = j^*\widetilde{T}j_*j^*T.
\]
The functors $i^!$ and $j_*$ have right-adjoint functors given by
\[
  i_? = Ti_*\widetilde{S} = Ti_*i^*\widetilde{T}i_*
  \;\;\;\; \mbox{and} \;\;\;\;
  j^! = Uj^*\widetilde{T} = j^*Tj_!j^*\widetilde{T}.
\]
\end{Proposition}

\begin{proof}
This can be proved directly, for instance
\renewcommand{\arraystretch}{1.2}
\[
  \begin{array}{rcll}
    \sT(i_!X,Y)
    & =     & \sT(\widetilde{T}i_*SX,Y)
            & \mbox{by definition} \\
    & \cong & \sT(Y,i_*SX)^{\vee}
            & \mbox{$\widetilde{T}$ left Serre functor} \\
    & \cong & \sS(i^*Y,SX)^{\vee}
            & \mbox{$i^*$ left-adjoint of $i_*$} \\
    & \cong & \sS(X,i^*Y),
            & \mbox{$S$ right Serre functor},
  \end{array}
\]
\renewcommand{\arraystretch}{1}
$\!\!$and similarly for the other cases.
\end{proof}

This permits the proof of the main result of this note.

\begin{Theorem}
\label{thm:main}
There are recollements
\[
\xymatrix
{
&&\sU
    \ar[rr]^{j_!} & & 
  \sT 
    \ar[rr]^{i^*}
    \ar@/_1.5pc/[ll]_{j^?} \ar@/^1.5pc/[ll]^{j^*} & &
  \sS
    \ar@/_1.5pc/[ll]_{i_!} \ar@/^1.5pc/[ll]^{i_*} & &
}
\]
and
\[
\xymatrix
{
&&\sU
    \ar[rr]^{j_*} & & 
  \sT 
    \ar[rr]^{i^!}
    \ar@/_1.5pc/[ll]_{j^*} \ar@/^1.5pc/[ll]^{j^!} & &
  \sS \lefteqn{.}
    \ar@/_1.5pc/[ll]_{i_*} \ar@/^1.5pc/[ll]^{i_?} & &
}
\]
\end{Theorem}

\begin{proof}
Proposition \ref{pro:adjoints} implies that there is
\[
\xymatrix
{
&&
  \sT 
    \ar[rr]^{i^*} & &
  \sS
    \ar@/_1.5pc/[ll]_{i_!} \ar@/^1.5pc/[ll]^{i_*} & &
}
\]
where $(i_!,i^*)$ and $(i^*,i_*)$ are pairs of adjoint functors.  The
functor $i_*$ is fully faithul, and it follows from \cite[prop.\
2.7]{Miyachi} or \cite[prop.\ 1.14]{Heider} that there is a
recollement
\[
\xymatrix
{
&&\Ker i^*
    \ar@{^{(}->}[rr] & & 
  \sT 
    \ar[rr]^{i^*}
    \ar@/_1.5pc/[ll] \ar@/^1.5pc/[ll] & &
  \sS.
    \ar@/_1.5pc/[ll]_{i_!} \ar@/^1.5pc/[ll]^{i_*} & &
}
\]
It is standard recollement theory that $\Ker i^* = \EssIm j_!$, see
\cite[thm.\ 1]{Konig} or \cite[rmk.\ 1.5(8)]{Heider}, and $j_!$ can be
used to replace $\EssIm j_!$ with $\sU$, so the first recollement of
the theorem,
\[
\xymatrix
{
&&\sU
    \ar[rr]^{j_!} & & 
  \sT 
    \ar[rr]^{i^*}
    \ar@/_1.5pc/[ll]_{j^?} \ar@/^1.5pc/[ll]^{j^*} & &
  \sS,
    \ar@/_1.5pc/[ll]_{i_!} \ar@/^1.5pc/[ll]^{i_*} & &
}
\]
follows.  The functors from $\sT$ to $\sU$ must be $j^?$ and $j^*$
since, by the definition of recollements, they are the left- and the
right-adjoint of the functor $j_!$ from $\sU$ to $\sT$.

The second recollement of the theorem can be obtained by the dual
procedure.
\end{proof}


\end{document}